\documentclass[12pt]{article}
\usepackage[english,british]{babel}
\usepackage[active]{srcltx}

\usepackage{amsfonts,amssymb,amstext}
\usepackage{amsmath}
\usepackage{amstext}
\usepackage{amsthm}

\makeatletter

\newcommand{\tr}{\hbox{ tr}}

\begin{document}
\centerline{\Large \bf On Riemann's mapping theorem.}
\vskip0.5cm
\centerline{\bf Ashot\,Vagharshakyan}
\vskip0.5cm {\small Institute
of Mathematics Armenian National Academy of
   Sciences, Bagramian 24-b, Yerevan, Armenia.}

{\small e-mail: vagharshakyan@yahoo.com} \vskip0.5cm

{\bf Abstract.}
    In this paper we give a new proof of Riemann's well known mapping
theorem. The suggested method permits to prove an analog of that
theorem for three dimensional case.


\section{Introduction}

    By Liouville's theorem, see [2], p. 130, in three dimensional case,
only superposition of isometric, dilatation and inverse
transformations are conformal. To get an analog for Riemann's
mapping theorem, one introduce a family of mappings named quasi -
conformal. This family is wider, nevertheless we do not have a
natural analog of conformal mappings like of two dimensional case.
In this paper we introduce a new family of mappings, named weak -
conformal and obtain more natural generalization of Riemann's
theorem.

    The proof of the main result of the present paper is interesting for two
dimensional case too. Actually, we give a new proof of Riemann's
classical theorem, where the specific properties of complex analysis
do not used. This permits us to prove a similar theorem in three
dimensional case.

\section{Classes of mappings}

    For any matrix $M=\{a_{ij}\}$ with eigenvalues $\lambda_1,\quad \dots,\lambda_n$ denote by
\[
|M|^2=\sum_{i=1}^n\sum_{j=1}^na_{ij}^2=\sum_{i=1}^n|\lambda_i|^2,
\]
\[
\tr(M)=\sum_{i=1}^na_{ii}=\sum_{i=1}^n\lambda_i
\]
and \[\det(M)=\prod_{i=1}^n\lambda_i.\] Let $\varphi(x,y,z)=(A, B,
C)$ be a continuously differentiable mapping. Denote by
\[
J=\left(\begin{array}{cccc}
A'_x & A'_y & A'_z\\
B'_x & B'_y & B'_z\\
C'_x & C'_y & C'_z\\
\end{array}\right)
\]
the Jacobi matrix. Let $G=J^*J$. We have
\[
|\varphi(\vec{x}+\Delta\vec{x})-\varphi(\vec{x})|^2=\left(J\Delta
\vec{x},J\Delta \vec{x}\right)+o(|\Delta \vec{x}|^2)= \left(\Delta
\vec{x},G\Delta \vec{x}\right)+o(|\Delta \vec{x}|^2).
\]

\vskip0.2cm {\bf Definition 1.} A continuously differentiable one to
one mapping
\[
\varphi:\Omega_1\rightarrow \Omega_2
\]
of the domain $\Omega_1\subset R^3$ on $\Omega_2\subset R^3$ is
conformal if for each point $\vec{x}\in\Omega_1$ there is a number
$M(\vec{x})$ such that
\[
|\varphi(\vec{x}+\Delta\vec{x})-\varphi(\vec{x})|=
M(\vec{x})\left|\Delta \vec{x}\right|+o(|\Delta \vec{x}|).
\]

 \vskip 0.2cm {\bf Lemma 1.} Let $\varphi$ be a continuously differentiable mapping with the
Jacobi matrix $J$ and $G=J^*J$. Then $\varphi$ is conformal if and
only if
\[
27\det(G)=\tr^3\left(G\right).
\]
\vskip0.2cm {\bf Proof.} The eigenvalues $\lambda_1, \lambda_2,
\lambda_3$ of the matrix $G$ are nonnegative. The lemma's condition
means that
\[
\left(\frac{\lambda_1+\lambda_2+\lambda_3}{3}\right)^3=\lambda_1\lambda_2\lambda_3
\]
This equality is valid only if all eigenvalues are the same, i. e.
$\lambda_1=\lambda_2=\lambda_3$. \vskip0.2cm
    {\bf Example.} Let us consider the inverse
transformation, which the point $(x,y,z)\neq(0,0,0)$ maps to
$(A,B,C)$, where
\[
(A, B, C)=\left(\frac{x}{x^2+y^2+z^2},
\frac{y}{x^2+y^2+z^2},\frac{z}{x^2+y^2+z^2},\right)
\]
We have
\[
J=\left(\begin{array}{cccc}
A'_x & A'_y & A'_z\\
B'_x & B'_y & B'_z\\
C'_x & C'_y & C'_z\\
\end{array}\right)=
\]
\[
=\frac{1}{(x^2+y^2+z^2)^2}\left(\begin{array}{cccc}
-x^2+y^2+z^2 & -2xy & -2xz\\
-2xy & x^2-y^2+z^2 & -2yz\\
-2xz & -2yz & x^2+y^2-z^2\\
\end{array}\right)
\]
Consequently,
\[
G=\frac{1}{(x^2+y^2+z^2)^2}\left(\begin{array}{cccc}
1 & 0 & 0\\
0 & 1 & 0\\
0 & 0 & 1\\
\end{array}\right)
\]
So, the condition of lemma 1 is satisfy and hence this mapping is
conformal. \vskip0.2cm {\bf Definition 2.} A quasi-conformal mapping
is a continuously differentiable homeomorphism
\[
\varphi:\Omega_1\rightarrow \Omega_2
\]
for which the ball
\[
B(\vec{x},r)=\{\vec{y};\,|\vec{x}-\vec{y}|<r\}
\]
maps to
\[
\{\varphi(\vec{y});\,\vec{y}\in
B(\vec{x},r)\}=\{\varphi(\vec{x})+J(\vec{y}-\vec{x}); \,\vec{y}\in
B(\vec{x},r)\}+o(r^2)
\]
and the ratio of the main diagonals  of the ellipsoid
\[
\{J(\vec{z});\,|\vec{z}|=r)\}
\]
are uniformly bounded for all points $\vec{x}\in \Omega$.
\vskip0.2cm

    In this paper we introduce a new family of mappings, which are generalizations of conformal mappings.
For those mappings, which we name weak - conformal, we have an
analog of Riemann's mapping theorem. \vskip0.2cm {\bf Definition 3.}
A weak-conformal mapping is a continuously differentiable
homomorphism
\[
\varphi:\Omega_1\rightarrow \Omega_2
\]
for which the ball
\[
B(\vec{x},r)=\{\vec{y};\,|\vec{x}-\vec{y}|<r\}
\]
maps to
\[
\{\varphi(\vec{y});\,\vec{y}\in
B(\vec{x},r)\}=\{\varphi(\vec{x})+J(\vec{y}-\vec{x}); \,\vec{y}\in
B(\vec{x},r)\}+o(r^2)
\]
and the the main diagonals  of the ellipsoid
\[
\{J(\vec{z});\,|\vec{z}|=r)\}
\]
form geometric progression for all points $\vec{x}\in \Omega$.
\vskip 0.2cm {\bf Lemma 2.} Let $\varphi$ be a continuously
differentiable mapping with Jacobi matrix $J$. Then it is weak -
conformal if and only if
\[
\left(\tr^2(G)-|G|^2\right)^3=8\det(G)\tr^3(G),
\]
where $G=J^*J$.
\vskip0.2cm {\bf Proof.} In terms of eigenvalues of
the matrix $G=J^*J$, this condition one can write as follows
\[
\lambda_1\lambda_2\lambda_3\left(\lambda_1+\lambda_2+\lambda_3\right)^3=
\left(\lambda_1\lambda_2+\lambda_3\lambda_1+\lambda_3\lambda_2\right)^3.
\]
So,
\[
(\lambda_1\lambda_2\lambda_3-\lambda_3^3)\left(\lambda_1+\lambda_2+\lambda_3\right)^3-
\left(\lambda_1\lambda_2+\lambda_3\lambda_1+\lambda_3\lambda_2\right)^3+
\left(\lambda_1\lambda_3+\lambda_2\lambda_3+\lambda_3^2\right)^3=0.
\]
After simple transformations we get
\[
(\lambda_1\lambda_2-\lambda_3^2)\left(\left(\lambda_1\lambda_2+\lambda_3\lambda_1+\lambda_3\lambda_2\right)
\left(\lambda_1^2+\lambda_2^2+\lambda_1\lambda_2\right)-
\lambda_1\lambda_2\left(\lambda_1+\lambda_2+\lambda_3\right)^2\right)=0.
\]
The last condition is equivalent to the following one
\[
\left(\lambda_1\lambda_2-\lambda_3^2\right)
\left(\lambda_1\lambda_3-\lambda_2^2\right)
\left(\lambda_3\lambda_2-\lambda_1^2\right)=0.
\]
Consequently, our condition means that eigenvalues of the matrix
$J^*J$ form a geometric progression.

\section{Green's function in $R^2$}

    In this section we introduce Green function and prove some of its properties.
\vskip0.2cm{\bf Definition 4.} Let $\Omega$ be a domain in $R^2$. A
function $G(\vec{x},\vec{y}),\quad \vec{x}\neq \vec{y}\in \Omega$ is
called Green function for the domain $\Omega$, if it satisfies the
following conditions:

1. $G(\vec{x},\vec{y})$ is continuous from below and
\[
G(\vec{x},\vec{y})>0,\quad \vec{x}\neq \vec{y} \in \Omega;
\]

2. for each fixed point $\vec{y}\in \Omega$ there is a harmonic
function $h(\vec{x},\vec{y}),\quad \vec{x}\in \Omega$ such that
\[
G(\vec{x},\vec{y})=\frac{1}{2\pi}\log\frac{1}{|\vec{x}-\vec{y}|}+h(\vec{x},\vec{y});
\]

3. if $u(\vec{x})$ is an arbitrary harmonic function defined in
$\Omega$ and satisfying the condition
\[
u(\vec{x})\leq G(\vec{x},\vec{y}),\quad \vec{x}\in \Omega\setminus
\{\vec{y}\},
\]
then
\[
u(\vec{x})\leq 0,\quad \vec{x}\in \Omega.
\]
\vskip0.2cm
    In three dimension case, if $\Omega \subset R^3$, the Green's function can be defined by the same way
replacing the second condition to the following one
\[
G(\vec{x},\vec{y})=\frac{1}{4\pi |\vec{x}-\vec{y}|}
+h(\vec{x},\vec{y}),\quad \vec{x}\in \Omega\setminus \{\vec{y}\}.
\]

    It is well known, that if the boundary of a domain $\Omega \neq R^2$ has
positive capacity, then it has unique Green function, see [4], p.
138. In particular, if $\Omega \neq R^2$ is simply connected then it
has Green function.

    For a fixed point $\vec{y}\in \Omega$ and for an arbitrary number $0<t<+\infty$ let us
denote
\[
\Omega_t=\{\vec{x};\quad G(\vec{x},\vec{y})>t\}.
\]
Note that for arbitrary value of $t>0$ the domain $\Omega_t$ is
connected and its Green function is $G(\vec{x},\vec{y})-t$.
\vskip0.5cm {\bf Lemma 3.} Let $u(\vec{x}),\quad \vec{x}\in \Omega,$
be a harmonic function and
$$\{\vec{x};\quad |\vec{x}-\vec{x}_0|\leq r\}\subset \Omega.$$ If for some point
$\vec{x},\quad |\vec{x}-\vec{x}_0|=r,$ we have
\[
u(\vec{x})= \inf\{u(\vec{y});\quad |\vec{y}-\vec{x}_0|<r\}
\]
then
\[
|\nabla u(\vec{x})|\geq \frac{1}{2r}(u(\vec{x}_0)-u(\vec{x})).
\]
\vskip0.2cm {\bf Proof.} For arbitrary $0\leq \varphi<2\pi$ and
$0<t<1$ we have
\[
\frac{r^2-(rt)^2}{|re^{i\varphi}-t(\vec{x}-\vec{x}_0)|^2}\geq
\frac{r^2-(rt)^2}{(r+rt)^2}\geq\frac{1-t}{2}.
\]
So,
\[
\frac{u(\vec{x}_0+t(\vec{x}-\vec{x}_0))-u(\vec{x})}{|\vec{x}_0+t(\vec{x}-\vec{x}_0)-\vec{x}|}=
\]
\[
=\frac{1}{2\pi(r-rt)}\int_0^{2\pi}\frac{r^2-(rt)^2}{|re^{i\varphi}-
t(\vec{x}-\vec{x}_0)|^2}(u(\vec{x}_0+re^{i\varphi})-u(\vec{x}))d\varphi\geq
\]
\[
\geq \frac{1}{4\pi
r}\int_0^{2\pi}(u(\vec{x}_0+re^{i\varphi})-u(\vec{x}))d\varphi=
\frac{1}{2r}(u(\vec{x}_0)-u(\vec{x})).
\]
Passing to the limit if $t\to 1-0$ we get the required result.
\vskip0.2cm
    {\bf Remark.} The analogous result is true in $R^3$.

\vskip0.2cm {\bf Theorem 1.} Let $G(\vec{x},\vec{y})$ be a Green
function for a simply connected domain $\Omega$ in $R^2$. Then
\[
|\nabla G(\vec{x},\vec{y})|\neq 0, \quad \vec{x}\in \Omega\setminus
\{\vec{y}\}.
\]
\vskip0.2cm {\bf Proof.} Let us assume, that at some point
$\vec{x}_0\in \Omega$ we have
$$
|\nabla G(\vec{x}_0,\vec{y})|=0.
$$
Denote
$$
\Omega^+=\{\vec{x},\quad G(\vec{x},\vec{y})>G(\vec{x}_0,\vec{y})\}
$$
and
$$
\Omega^-=\{\vec{x},\quad G(\vec{x},\vec{y})<G(\vec{x}_0,\vec{y})\}.
$$
Note that the domain $\Omega^+$ is connected. For sufficiently small
number $\epsilon>0,$ such that the following condition
$$
B(\vec{x}_0,\epsilon)\subset\Omega\setminus \{\vec{y}\},
$$
holds, we consider the open set
\[
A=B(\vec{x}_0,\epsilon)\setminus \{\vec{x};\quad
G(\vec{x},\vec{y})=G(\vec{x}_0,\vec{y})\}.
\]

    The set $A$ consists of an even number of components. Otherwise, we
could find a point $$\vec{x}_1\in
B(\vec{x}_0,\epsilon)\cap\{\vec{x};\quad
G(\vec{x},\vec{y})=G(\vec{x}_0,\vec{y})\}$$ in some neighborhood of
which the function $G(\vec{x},\vec{y})-G(\vec{x}_0,\vec{y})$ would
preserve its sign. So, $\vec{x}_1$ would be the point of local
extremum, which is impossible for the nonconstant harmonic function
$G(\vec{x},\vec{y}),\quad \vec{x}\in B(\vec{x}_0,\epsilon)$.
Moreover, the set $A$ cannot have only two components. Indeed if it
had two components then the boundary $\partial \Omega^+ \cap
B(\vec{x}_0,\epsilon)$ would be smooth. Consequently, by lemma 1 we
would have
$$
|\nabla G(\vec{x}_0,\vec{y})|>0.
$$
Thus, the domain $A$ has at least four connected components. This
implies that the open set $\Omega^-$ consists of more than two
connected components. Since our domain is simply connected, one of
those components has the boundary, completely laying inside of
$\partial\Omega^+$. On that connected component the function
$G(\vec{x},\vec{y})$ is identically constant end equal
$G(\vec{x}_0,\vec{y})$. This is a contradiction.

    Note that in theorem 1, the condition "simply connected", is essential.
Indeed, for the domain $\{\vec{x};\quad 1<|\vec{x}|<2\}$ theorem 1
does not valid.

\section{New proof of Riemann's theorem}

    In this section we give a new proof of Riemann's well known theorem on conformal
mapping. In this proof we do not use methods of complex analysis.

    Let us denote
\[
D(\vec{y},r)=\{\vec{x};\quad |\vec{x}-\vec{y}|<r\}.
\]
\vskip0.5cm {\bf Theorem 2.} Let $\Omega$ be a simply connected
domain in $R^2$. If $\Omega\neq R^2$ then there is a one to one
conformal mapping
\[
\varphi: \Omega \rightarrow D
\]
of the domain $\Omega$ on the unit disk $D=D(\vec{0},1)$.

\vskip0.5cm {\bf Proof.} Let us fix a point $\vec{y}\in \Omega$ and
$G(\vec{x},\vec{y})$ be Green function of the domain $\Omega$. Let
us consider the following dynamical system in $\Omega\setminus
\{\vec{y}\}$
\[
\frac{d\vec{x}(t)}{dt}=-\frac{\nabla
G(\vec{x}(t),\vec{y})}{2\pi|\nabla G(\vec{x}(t),\vec{y})|^2}e^{2\pi
G(\vec{x}(t),\vec{y})},\quad 0<t<1.\quad\quad (1)
\]

    For an arbitrary solution of this equation we have
\[
\frac{d}{dt}\left(e^{-2\pi G(\vec{x}(t),\vec{y})}\right)=-2\pi
e^{-2\pi G(\vec{x}(t),\vec{y})}\left(\nabla G(\vec{x}(t),\vec{y}),
\frac{d\vec{x}(t)}{dt}\right)=1.
\]
Consequently,
\[
G(\vec{x}(t),\vec{y})=\frac{1}{2\pi}\ln\frac{1}{t},\quad 0<t<1.
\]

    In the neighborhood of each point $\vec{x}\in \Omega\setminus \{\vec{y}\}$
the equation (1) has a unique solution passing through the point
$\vec{x}$, see [1] p. 19.

    In the neighborhood of the point $\vec{y}$ the equation (1) may be
written in the following form
\[
\frac{d\vec{x}(t)}{dt}=\frac{\vec{x}(t)-\vec{y}}{|\vec{x}(t)-\vec{y}|}\exp\{2\pi
h(\vec{y},\vec{y})\}+o(t),\quad t\to 0.
\]
So, for each solution of our equation we have
\[
\vec{x}(t)=\vec{y}+\vec{a}t\exp\{2\pi
h(\vec{y},\vec{y})\}+o(t),\quad t\to 0,
\]
where $\vec{a}$ is a vector with norm one.

    Consequently, for each point $x\in \Omega\setminus \{\vec{y}\}$ we can
find a unique vector $\vec{a}=\vec{a}(\vec{x})$, such that there is
a solution $\vec{x}(t)$ of our equation which passes through the
point $\vec{x}$ and at the same time in the neighborhood of the
point $\vec{y}$ satisfies the condition
\[
\lim_{t\to 0}\frac{\vec{x}(t)-\vec{y}}{t}=\vec{a}\exp\{2\pi
h(\vec{y},\vec{y})\}.
\]

    Let us define the mapping
$$
\varphi : \Omega \rightarrow D
$$
as follows, $\varphi(\vec{y})=0$ and for the arbitrary point
$\vec{x}\in \Omega\setminus \{\vec{y}\}$ we put
$$
\varphi(\vec{x})=\vec{a}(\vec{x})e^{-2\pi G(\vec{x},\vec{y})}.
$$
It is obvious, that $\varphi(\vec{x})$ is a one to one mapping and
$\varphi(\Omega)=D$.

    Recall some facts about the constructed mapping, which
permit to assert that it is conformal.

    Let us take two solutions
$$\vec{x}(t),\quad \vec{x}_1(t)$$
of the equation (1). We denote by $\alpha$ the angle between the
vectors $\vec{a}(\vec{x}(t))$ and $\vec{a}(\vec{x}_1(t))$. For
arbitrary numbers $0<t_0<t_1<1$ denote by $U$ the domain bounded by
the curves
$$\gamma_1=\{\vec{x}(t);\quad t_0<t<t_1\},\quad \gamma_2=\{\vec{x}_1(t);\quad
t_0<t<t_1\}$$ and
$$\gamma_3=\{\vec{x};\quad G(\vec{x},\vec{y})=G(\vec{x(t_0)},\vec{y})\},
\quad \gamma_4=\{\vec{x};\quad
G(\vec{x},\vec{y})=G(\vec{x(t_1)},\vec{y})\}.$$ Let
$\vec{m}(\vec{x})$ be the unite outer normal to the boundary of the
domain $U$ at the point $\vec{x}\in \partial U$. For an arbitrary
point $\vec{x}\in \gamma_1\cup\gamma_2$ we have
\[
\left(\frac{d\vec{x}(t)}{dt}, \vec{m}(\vec{x}(t))\right)=0,\quad
t_0<t<t_1.
\]
Consequently,
$$
\left(\nabla G(\vec{x},\vec{y}),
\vec{m}(\vec{x})\right)=\frac{\partial G(\vec{x},\vec{y})}{\partial
\vec{m}}=0.
$$
If $\vec{x}\in \gamma_3$ then we have
$\vec{m}(\vec{x})=-\vec{n}(\vec{x})$, where $\vec{n}(\vec{x})$ is
the outer normal to the boundary of the domain $\{\vec{x};\quad
G(\vec{x},\vec{y})>t_0\}.$ If $\vec{x}\in \gamma_4$ then we have
$\vec{m}(\vec{x})=\vec{n}(\vec{x})$, where $\vec{n}(\vec{x})$ is the
outer normal to the domain $\{\vec{x};\quad G(\vec{x},\vec{y})>t\}.$
Therefore,
$$
\int_{\gamma_3}\frac{\partial G(\vec{x},\vec{y})}{\partial
\vec{n}}ds=\int_{\gamma_4}\frac{\partial
G(\vec{x},\vec{y})}{\partial \vec{n}}ds.
$$
Passing to the limit we get
$$
\alpha=2\pi\lim_{t_0\to +0}\int_{\gamma_3}\frac{\partial
G(\vec{x},\vec{y})}{\partial
\vec{n}}ds=2\pi\int_{\gamma_4}\frac{\partial
G(\vec{x},\vec{y})}{\partial \vec{n}}ds.
$$

    From definition of the mapping $\varphi(\vec{x})$ we have
$$
|\varphi(\vec{x}(t))-\varphi(\vec{x}_1(t)|=|t(\vec{x})||\vec{a}(\vec{x}(t))-\vec{a}(\vec{x}_1(t))|=
$$
$$
=|t(\vec{x})|\left|2\pi\int_{\gamma_4}\frac{\partial
G(\vec{x},\vec{y})}{\partial \vec{n}}ds\right|=
$$
$$
= 2\pi\frac{\partial G(\vec{x}(t),\vec{y})}{\partial
\vec{n}}\exp\{-2\pi
G(\vec{x}(t),\vec{y})\}|\vec{x}(t)-\vec{x}_1(t)|+o(|\vec{x}(t)-\vec{x}_1(t)|).
$$
Further on we can write
\[
|\varphi(\vec{x}(t+\Delta t))-\varphi(\vec{x}(t))|+o(|\Delta
t|)=|\Delta t|+o(|\Delta t|)=
\]
\[
=|\vec{x}(t+\Delta
t)-\vec{x}(t)|\left|\frac{d\vec{x}(t)}{dt}\right|^{-1}=
\]
\[
=2\pi\frac{\partial G(\vec{x}(t),\vec{y})}{\partial
\vec{n}}\exp\{-2\pi G(\vec{x}(t),\vec{y})\}|\vec{x}(t+\Delta
t)-\vec{x}(t)|+o(|\Delta t|).
\]

    Let $r>0$ be sufficiently small. We choose $\Delta t$ and $\vec{x}_1(t)$ such
that the equalities
\[
|\vec{x}(t+\Delta t)-\vec{x}(t)|=|\vec{x}(t)-\vec{x}_1(t)|=r
\]
hold. The vectors $$\vec{x}(t+\Delta t)-\vec{x}(t)$$ and
$$\vec{x}(t)-\vec{x}_1(t)$$ are orthogonal. Consequently, the image of the
disk $D(\vec{x}(t),r)$ is a circle, as a first approximation, once
if the orthogonal vectors
$$\varphi(\vec{x}(t+\Delta t))-\varphi(\vec{x}(t))$$ and
$$\varphi(\vec{x}(t))-\varphi(\vec{x}_1(t))$$ satisfy the condition
$$
|\varphi(\vec{x}(t+\Delta t))-\varphi(\vec{x}(t))|=
|\varphi(\vec{x}(t))-\varphi(\vec{x}_1(t))|+o(r).
$$
The last condition holds since
$$
|\varphi(\vec{x}(t+\Delta t))-\varphi(\vec{x}(t))|=
2\pi\frac{\partial G(\vec{x}(t),\vec{y})}{\partial
\vec{n}}\exp\{-2\pi G(\vec{x}(t),\vec{y})\}r+o(r)
$$
and
$$
|\varphi(\vec{x}(t))-\varphi(\vec{x}_1(t)|= 2\pi\frac{\partial
G(\vec{x}(t),\vec{y})}{\partial \vec{n}}\exp\{-2\pi
G(\vec{x}(t),\vec{y})\}r+o(r).
$$
\vskip0.2cm
    {\bf Remark.} For constructed mapping at the points $\vec{x}\in\Omega\setminus\{\vec{y}\}$ we have
$$
|\varphi'(\vec{x})|=2\pi|\nabla G(\vec{x},\vec{y})|\exp\{-2\pi
G(\vec{x},\vec{y})\}.
$$
At the point $\vec{y}$ we have
$$
|\varphi'(\vec{y})|=2\pi\exp\{-2\pi h(\vec{y},\vec{y})\}.
$$

\section {Green's function in $\bf R^3$}

\vskip0.2cm {\bf Definition 5.} We say that a domain $\Omega\subset
R^3$ is simply connected if

    1. for an arbitrary bounded domain $\Omega_1\subset
R^3$ if we have $\partial\Omega_1\subset \Omega$ then it follows
$\Omega_1\subset \Omega$;

    2. an arbitrary closed curve laying in domain $\Omega$ permits
continuous deformation in domain $\Omega$ to the point.

\vskip0.2cm {\bf Lemma 4.} Let $\Omega$ be a simply connected domain
in $R^3$. Let $\Omega$ be a bounded domain with smooth boundary and
$G(\vec{x},\vec{y})$ is its Green function. Then
\[
\nabla G(\vec{x},\vec{y})\neq 0,\quad \vec{x}\in \Omega\setminus
\{\vec{y}\}.
\]

\vskip0.2cm {\bf Proof}. Since the boundary of our domain is smooth
so, we have
\[
\nabla G(\vec{x},\vec{y})\neq 0,\quad \vec{x}\in
\partial\Omega.
\]

    Let us assume that
\[
\{\vec{x};\,\vec{x}\in \Omega, \,\,\nabla
G(\vec{x},\vec{y})=0\}\neq\emptyset.
\]
    Let $0<t_0<\infty$ be the biggest number for which there is a point $x_0\in \Omega$ such
that $G(\vec{x}_0,\vec{y})=t_0$ and
\[
\nabla G(\vec{x}_0,\vec{y})=0.
\]
Denote
$$
\Omega^+=\{\vec{x},\quad G(\vec{x},\vec{y})>t_0\}
$$
and
$$
\Omega^-=\{\vec{x},\quad G(\vec{x},\vec{y})<t_0\}.
$$
Note that the domain $\Omega^+$ is connected. If $\Omega^-$ does not
connected we come to the contradiction like of two dimensional case.

    It turns out, that in three dimensional case, it is possible that the domain $\Omega^-$
is connected too.

    In the domain
\[
\{\vec{x}; \vec{x}\in \Omega,\,G(\vec{x},\vec{y})>t_0\}
\]
we consider the following dynamic system
\[
\frac{d\vec{x}(t)}{dt}=-\frac{\nabla
G(\vec{x}(t),\vec{y})}{4\pi|\nabla
G(\vec{x}(t),\vec{y})|^2}G^{2}(\vec{x}(t),\vec{y}),\quad t_0<t.
\]

    For an arbitrary solution of this equation we have
\[
\frac{d}{dt}\frac{1}{G(\vec{x}(t),\vec{y})}=-\frac{1}{G^2(\vec{x}(t),\vec{y})}\left(\nabla
G(\vec{x}(t),\vec{y}),\,\frac{d\vec{x}(t)}{dt}\right)=\frac{1}{4\pi}.
\]
Consequently, we have
\[
G(\vec{x}(t),\vec{y})=\frac{1}{4\pi t},\quad 0<t<\infty.
\]

So, for each $0<\epsilon$, the solutions of this equation generate
the following transformation
\[
x(\infty)\rightarrow x(t_0+\varepsilon)
\]
which settle a one to one correspondence between the points of the
manifolds $\partial \Omega$ and $\partial \Omega_{t_0+\epsilon}$.
Hence, those manifolds are homotopic equivalent.

    In the domain
\[
\{\vec{x},\quad G(\vec{x},\vec{y})>t_0+\varepsilon\}
\]
there is a smooth closed curve $\gamma_1$, which passes through the
points $x_0$ and $y$.

    For sufficiently small $\epsilon>0$ the plane orthogonal to the
curve $\gamma_1$ at the point $x_0$, cut a closed curve $\gamma_2$
on the boundary $\partial\Omega_{t_0+\epsilon}$ which have nonzero
index in compare to the curve $\gamma_1$.

    Since $\partial \Omega_{t_0+\epsilon}$ and $\partial \Omega$ are
homotopic equivalent so, the curve $\gamma_2$, by continuously
deformation, staying on the boundary $\partial
\Omega_{t_0-\epsilon}$, is possible to tie up to a point.

    This is a contradiction since each curve on the boundary
$\partial \Omega_{t_0+\epsilon}$ having sufficiently small diameter,
has zero index in compare with the curve $\gamma_1$.

\section {Weak-conformal mapping in $\bf R^3$}

    In this section we prove the main result of the paper.
\vskip0.2cm {\bf Theorem 3.} Let $\Omega$ be a simply connected
domain in $R^3$. If $\Omega$ is a bounded and has smooth boundary
then there is a one to one weak - conformal mapping
\[
\varphi: \Omega \rightarrow B
\]
of the domain $\Omega$ onto the unit ball $B=\{x\in R^3;\quad
|x|<1\}$.

\vskip0.2cm {\bf Proof}. We consider the following dynamic system
\[
\frac{d\vec{x}(t)}{dt}=-\frac{\nabla
G(\vec{x}(t),\vec{y})}{4\pi|\nabla
G(\vec{x}(t),\vec{y})|^2}G^{2}(\vec{x}(t),\vec{y}),\quad 0<t<\infty.
\]

    In neighborhood of the point $\vec{y}$ we have
\[
-\frac{\nabla G(\vec{x},\vec{y})}{4\pi|\nabla
G(\vec{x},\vec{y})|^2}G^{2}(\vec{x},\vec{y})=
\]
\[
=\left(\frac{\vec{x}-\vec{y}}{4\pi |\vec{x}-\vec{y}|^3}-\nabla
h\right)\left(\frac{1}{4\pi |\vec{x}-\vec{y}|}+h\right)^2
\]
\[
\left(\frac{1}{16\pi^2|\vec{x}-\vec{y}|^4}-\frac{(\vec{x}-\vec{y},\,\nabla
h)}{2\pi |\vec{x}-\vec{y}|^3}+|\nabla h|^2\right)^{-1}=
\]
\[
=\left(\frac{\vec{x}-\vec{y}}{|\vec{x}-\vec{y}|}-4\pi|\vec{x}-\vec{y}|^2\nabla
h\right)\frac{\left(1+4\pi |\vec{x}-\vec{y}|h\right)^2}{
1-8\pi|\vec{x}-\vec{y}|(\vec{x}-\vec{y},\,\nabla
h)+16\pi^2|\vec{x}-\vec{y}|^4|\nabla h|^2}=
\]
\[
=\frac{\vec{x}-\vec{y}}{|\vec{x}-\vec{y}|}+8\pi
h(\vec{y},\vec{y})(\vec{x}-\vec{y})+O(|\vec{x}-\vec{y}|^2)
\]
So, for each solution of our equation we have
\[
\vec{x}(t)=\vec{y}+\vec{a}t+4\pi
\vec{a}t^2h(\vec{y},\vec{y})+o(t^2),\quad t\to 0,
\]
where $\vec{a}$ is a vector with norm one.

    Consequently, for each point $\vec{x}\in \Omega\setminus \{\vec{y}\}$ we can
find the unique vector $\vec{a}=\vec{a}(\vec{x})$ of unit norm, such
that a solution $\vec{x}(t)$ passes through the point $\vec{x}$ and
\[
\lim_{t\to 0}\frac{\vec{x}(t)-\vec{y}}{t}=\vec{a}.
\]
By definition the vector $\vec{a}(\vec{x}(t))$ is the same for all
values of $0<t<\infty$.

    Let $\vec{x}=\vec{x}(t_0)$. We denote by
\[
\gamma(\vec{x})=\{\vec{x}(t);\quad t_0\leq t<\infty\}
\]
the curve begins of the point $\vec{x}$ and goes to the boundary of
the domain $\Omega$.

    If for each point $\vec{x}\in \Omega$ the curve $\gamma(\vec{x})$
has a finite length, then we can define the mapping
$$
\varphi : \Omega \rightarrow B
$$
as follows, $\varphi(\vec{y})=0$ and for the point $\vec{x}\in
\Omega\setminus \{\vec{y}\}$ we put
$$
\varphi(\vec{x})=\vec{a}(\vec{x})\exp\left\{-\int_{\gamma(\vec{x})}\sqrt{4\pi
|\nabla G(\vec{z},\vec{y})|}ds(\vec{z})\right\}.
$$
It is obvious, that $\varphi$ is a one to one mapping onto the unit
ball $B$.

    Now let us consider the properties of the constructed mapping.

    For arbitrary nonzero vector $\vec{a}$ let us denote by
$D_{\alpha}(\vec{a})$ the round cone with bisector $\vec{a}$ and the
spherical sector
\[
\left\{\vec{y};\quad |\vec{y}|=1,\,\,\vec{y}\in
D_{\alpha}(\vec{a})\right\}
\]
has area equal $\alpha.$

    Let us fix a point $\vec{x}\in \Omega$ and a number $0<\alpha<4\pi$.
For arbitrary numbers $0<t_0<t_1<\infty$ denote by $U$ the following
domain
\[
U=\bigcup_{\vec{a}(\vec{x}(t_0))\in D_{\alpha}(\vec{a}(\vec{x}))}
\{\vec{x}(t);\quad t_0<t<t_1\}
\]
We denote
\[
W(t_0,\alpha)=\left\{\vec{x}(t_0);\quad \vec{a}(\vec{x}(t_0))\in
D_{\alpha}(\vec{a}(\vec{x}))\right\}
\]
and
\[
W(t_1,\alpha)=\left\{\vec{x}(t_1);\quad \vec{a}(\vec{x}(t_1))\in
D_{\alpha}(\vec{a}(\vec{x}))\right\}
\]
By Green's formula we have
\[
\int_{W(t_0,\alpha)}\frac{\partial G(\vec{z},\vec{y})}{\partial
\vec{n}}ds(\vec{z})= \int_{W(t_1,\alpha)} \frac{\partial
G(\vec{z},\vec{y})}{\partial \vec{n}}ds(\vec{z}).
\]
Passing to the limit if $t_0\to 0$ we get
\[
\int_{W(t_1,\alpha)}\frac{\partial G(\vec{z},\vec{y})}{\partial
\vec{n}}ds(\vec{z})=\lim_{t_0\to
+0}\int_{W(t_0,\alpha)}\frac{\partial G(\vec{z},\vec{y})}{\partial
\vec{n}}ds(\vec{z})=\frac{\alpha}{4\pi}.
\]
Consequently, for small $\alpha$ we have
\[
s(W(t_1,\alpha))\frac{\partial G(\vec{x}(t_1),\vec{y})}{\partial
\vec{n}}=\frac{\alpha}{4\pi}+o(\alpha).
\]
    The vector $\vec{x}(t_1+\Delta t)-\vec{x}(t_1)$ is orthogonal to
the surface $W(t_1,\alpha)$ and
\[
|\vec{x}(t_1+\Delta t)-\vec{x}(t_1)|=
\left|\frac{dx(t_1)}{dt}\right||\Delta t|=\frac{G^2(x(t_1),y)}{4\pi
|\nabla G(x(t_1),y)|}|\Delta t|+o(|\Delta t|).
\]
We choose the parameters $\alpha$ and $\Delta t$ such that
\[
s(W(t_1,\alpha))=\pi|\vec{x}(t_1+\Delta t)-\vec{x}(t_1)|^2.
\]
This means that we have
\[
\frac{G^4(x(t_1),y)}{16\pi |\nabla G(x(t_1),y)|^2}|\Delta
t|^2=\frac{\alpha}{4\pi |\nabla G(x(t_1),y)|}+o(\alpha).
\]

    The vector $\varphi(\vec{x}(t_1+\Delta
t)-\varphi(\vec{x}(t_1))$ is orthogonal to the surface
$\varphi(W(t_1,\alpha))$.

    Note that the image of the subset $W(t_1,\alpha)$ is a round
sector on the sphere with the center at the point $\vec{0}$ and with
the radius $|\varphi(\vec{x}(t_1))|$. So, we have
\[
s(\varphi(W(t_1,\alpha)))=\alpha|\varphi(\vec{x}(t_1))|^2.
\]
Since
\[
|\varphi(\vec{x}(t_1))|=\exp\left\{-\int_{\gamma(\vec{x}(t_1))}\sqrt{4\pi
|\nabla G(\vec{z},\vec{y})|}ds(\vec{z})\right\}=
\]
\[
=\exp\left\{-\int_{t_1}^{\infty}\frac{G^2(x(t),y)}{\sqrt{4\pi
|\nabla G(x(t),y)|}}dt\right\}
\]
So,
\[
|\varphi(\vec{x}(t_1+\Delta t))-\varphi(\vec{x}(t_1))|=
|\varphi(\vec{x}(t_1))|\frac{G^2(\vec{x}(t_1),\vec{y})}{\sqrt{4\pi
|\nabla G(\vec{x}(t_1),\vec{y})|}}|\Delta t|+o(|\Delta t|).
\]
Consequently, we have
\[
s(\varphi(W(t_1,\alpha)))=\alpha|\varphi(\vec{x}(t))|^2= \pi
|\varphi(\vec{x}(t+\Delta t))-\varphi(\vec{x}(t))|^2+o(|\Delta
t|^2).
\]

    This relation is equivalent to the weak - conformal condition at the point $\vec{x}\in
\Omega$ for constructed mapping.

{\bf Remark.} For any point $\vec{x}\neq \vec{y}$ we have
\[
\lim_{\Delta t\to 0}\frac{|\varphi(\vec{x}(t+\Delta
t))-\varphi(\vec{x}(t))|}{|\vec{x}(t+\Delta t)-\vec{x}(t)|}=
\]
\[
=\sqrt{4\pi |\nabla
G(\vec{x},\vec{y})|}\exp\left\{-\int_{\gamma(\vec{x})}\sqrt{4\pi
|\nabla G(\vec{z},\vec{y})|}ds(\vec{z})\right\}
\]
where $\vec{x}=\vec{x}(t)$. Hence, we have
\[
\varphi(\vec{x})=\varphi(\vec{x})-\varphi(\vec{y})=\vec{x}-\vec{y}+o(|\vec{x})-\vec{y}|).
\]

\vskip 0.5cm {\bf References}

1. Ph. Hartman, Ordinary differential equations, New York, London,
Sydney 1964.

2. B. A. Dubrovin, S. P. Novikov, A. T. Fomenko, Modern geometry,
Moscow 1986.

3. V. Arnold, A. Varchenko, S. Gusein - Zade, The singularities of
differentiable mappings, Moscow 1982

4. W. Hayman, P. Kennedy, Subharmonic functions, London, New York,
San Francisco, 1976.

5. T. Ransford, Potential theory in the complex plane, Cambridge
University press 1995.

6. L. Ahlfors, Mobius transformations in several dimensions,
University of Minnesota 1981.

\end{document}